\documentclass{article}
\usepackage{graphicx} 
\usepackage{amsmath} 
\usepackage{amssymb}
 
\newcommand{\la}{\lambda} 
\newcommand{\ga}{\gamma}
\newcommand{\al}{\alpha}

\newcommand{\De}{\Delta} 
\newcommand{\val}{\mathsf{val}}

\newcommand{\EE}{\mathbb{E}}

\newcommand{\RR}{\mathbb{R}}
\newcommand{\NN}{\mathbb{N}}
\newcommand{\QQ}{\mathbb{Q}}
\newtheorem{@theorem}{Theorem}

\newtheorem{Theorem}[@theorem]{Theorem}
\newtheorem{Lemma}[@theorem]{Lemma}
\newtheorem{lemma}[@theorem]{Lemma}

\newtheorem{Conjecture}{Conjecture}
\newtheorem{remarque}{Remark}

\title{Absorbing games with irrational values}
\author{Miquel Oliu-Barton\footnote{University Paris Dauphine - PSL, CNRS, CEREMADE, Paris, France.}}
\date{July 3, 2023}

\begin{document}

\maketitle

\begin{abstract}
Can an absorbing game with rational data have an irrational limit value?
Yes: In this note we provide the simplest examples where this phenomenon arises. That is, the following $3\times 3$ absorbing game \[
A = \begin{bmatrix}
1^* & 1^* & 2^* \\
1^* & 2^* & 0\phantom{^*} \\
2^* & 0\phantom{^*} & 1^*
\end{bmatrix},
\]
and a sequence of $2\times 2$ absorbing games whose limit values are $\sqrt{k}$, for all $k\in \NN^*$. 
Finally, 
we conjecture that any algebraic number can be represented as the limit value of an absorbing game. 
\end{abstract}

\section{Introduction}
In this note we pose a simple question: Can a rational absorbing game (i.e. one where payoffs and transitions are all rational) have an irrational limit value? The motivation for this question is twofold. First, as recently established in \cite{OB21}, for any rational stochastic game with $m$ actions per state and $K$ non-absorbing actions, the limit value is algebraic of degree $m^K$. Rational absorbing games could thus have an irrational limit value as soon as both players have $2$ actions per state. On the other hand, however, no rational absorbing game with an irrational limit value has ever been found. We fill this gap by providing the simplest possible examples of rational absorbing games with irrational values. These results contrast with Markov decision processes, turn-based stochastic games, and stochastic games in which transitions do not occur under optimal play, which all have rational limit values. Other classes of stochastic games are already known to have possibly irrational values (e.g. irreducible games \cite{CH14}).  



\paragraph{The Big Match.} 
Introduced by Everett in 1957 \cite{everett57}, the ``Big Match'' is the popular stochastic game: 
\[
\begin{array}{c|cc}
& \text{Player 2} \\
\text{Player 1} & L & R \\
\hline
T & 1^* & 0^* \\
B & 0\phantom{^*} & 1\phantom{^*} \\
\end{array}
\]
Like the Prisoner's dilemma for repeated games, the Big Match is definitely the most popular example of a stochastic game. Its success comes both from the simplicity of its representation (i.e., a real matrix and some $*$'s) and the complexity of its resolution by Blackwell and Ferguson \cite{BF68}. The game goes as follows: As long as Player 1 plays the bottom action `B', the stage payoffs are 0 or 1 depending on whether Player 2 is playing their left action `L' or right action `R', respectively. However, once Player 1 plays the top action `T', not only the current but also all future payoffs will be 1 or 0 depending, respectively, on whether Player 2 played `L' or `R' at that stage. 

\paragraph{Absorbing games.} Those are stochastic games in which the state can change at most once. 
Formally, an absorbing game is represented by three matrices $g,q,w\in \RR^{m\times n}$ where for each pair of actions $(i,j)$, $g_{ij}$ is for the non-absorbing payoff, $q_{ij}$ is the probability of absorption, and $w_{ij}$ is the absorbing payoff, which only matters if $q_{ij}>0$. An absorbing games with deterministic transitions is one that satisfies $q_{ij}\in \{0,1\}$ for all $(i,j)$. Provided that, in addition, $w_{ij}=g_{ij}$ for all $(i,j)$, like in the Big Match, the game can then be represented by a single matrix in which some payoffs are labelled with a $*$ to indicate that this stage payoffs is fixed from that stage onward. 

\paragraph{Literature.} The resolution of the Big Match motivated the formalization of absorbing games, due to Kohlberg \cite{kohlberg74}, who in addition proved the existence of the limit value, and provided a characterization using the derivative of the Shapley operator. Two additional characterizations for the limit value were obtained by Laraki \cite{laraki10} and Sorin and Vigeral \cite{SV13}. The complexity of the limit value was bounded in \cite{HKLMT11}, and improved in \cite{OB21}. General notes on stochastic games can be found in \cite{renaultnotes2, solan22,sorin03}. 





\section{Main result}

\begin{Theorem} A rational $m\times n$ absorbing game with deterministic transitions and $\min(m,n)<3$ has a rational limit value. 
\end{Theorem}

\begin{Theorem} There exists an integer $3\times 3$-matrix with stars with an irrational
limit value\footnote{The limit value of the game is the unique real solution of $P(z)=z^3-5z^2+10z-7$. This number is algebraic of degree $3$ as there is no smaller-degree polynomial to which it is a root.}, namely:  
\[
A = \begin{bmatrix}
1^* & 1^* & 2^* \\
1^* & 2^* & 0\phantom{^*} \\
2^* & 0\phantom{^*} & 1^*
\end{bmatrix}\,.
\]
\end{Theorem}
\begin{Theorem}\label{thm1} 
For any $k\in \NN^*$, the following absorbing game has limit value $\sqrt k$:
\[\begin{bmatrix}
0; (\frac{1}{k},k^*) & 1^*\\
1^* & k^*
\end{bmatrix},
\] where the entry $(g,(q,w^*))$ indicates the non-absorbing payoff $g$, the absorption probability $q$, and the absorbing payoff $w$. 
\end{Theorem}
\begin{remarque} Any algebraic number of degree $2$ is the limit value of $2\times 2$-absorbing game. Indeed, it is enough to consider affine transformations of the example of Theorem~\ref{thm1}. \end{remarque}


In view of our results, we further propose two open problems. 
\begin{Conjecture} For all $m\geq 1$, there exists a rational absorbing game with $m$ actions per player whose limit value is algebraic of order $m$. 
\end{Conjecture}
\begin{Conjecture} Any algebraic number can be represented as the limit value of a rational absorbing game. 
\end{Conjecture}
The first conjecture is equivalent to the bound in \cite{OB21} being tight in the class of absorbing games. This note solves the tightness of this bound only for the cases $m=2$ and $m=3$. 
The second is reminiscent of similar representation results, for example \cite{ASZ20} proved that any piece-wise rational fraction is the value of a polynomial matrix game, while \cite{VV16} proved that any compact semi-algebraic set is the projection of the set of Nash equilibrium payoffs of a game.  



\section{Proofs}
\subsection{Notation and useful known results} 
Throughout the paper we assume that state $1$ is the unique non-absorbing state, and that $v_\la$ and $v$ denote respectively the discounted value and their limit. For any stationary pair of strategies $(x,y)\in \De(I)\times \De(J)$, 
we denote by $\ga_\la(x,y)$ the expected payoff, starting from the initial state $1$. By \cite{shapley53}, we know that 
\begin{equation*}\label{noteq}
v_\la=\max_{x\in \De(I)}\min_{j\in J}\ga_\la(x,j)\,.
\end{equation*}
where $I=\{1,\dots,m\}$, $J=\{1,\dots,n\}$, and $\De(I)$ and $\De(J)$ denote, respectively, the set of probabilities over $I$ and $J$.
By finiteness, we also know that
$$v=\lim_{\la\to 0} \left(\max_{x\in \De(I)}\min_{j\in J}\ga_\la(x,j)\right)=
\max_{x\in \De(I)}\min_{j\in J}\left( \lim_{\la\to 0}\ga_\la(x,j)\right)\,.$$
Given an absorbing game $(g,q,w)\in \RR^{m\times n}$, 
define a parameterized matrix $W_\la(z)\in \RR^{m\times n}$ 
as follows: 
\begin{equation}\label{w}
W_\la(z)_{ij}=\la g_{ij} +(1-\la)q_{ij} w_{ij}-z(\la + (1-\la)q_{ij}),\quad \forall (i,j)\,.
\end{equation}

The following result will be used in the sequel. 
\begin{lemma}\label{lem} 
There exists 
a square sub-matrix $\dot{W}_\la(z)$ of $W_\la(z)$ so that 
\begin{itemize}
    \item $P(\la,z):=\det(\dot{W}_\la(z))$ is of degree at least $1$ in $z$ and satisfies $P(\la,v_\la)=0$; 
    \item The rows and columns of $\dot{W}_\la(z)$ correspond to the support of a pair of optimal stationary strategies. 
\end{itemize}
\end{lemma}
\textbf{Proof.} The matrix $W_\la(z)$ is, up to a strictly positive constant (i.e. $\la^K$, where $K$ is the number of non absorbing states) equal to the auxiliary matrix introduced in \cite{AOB19}. The two statements thus follow from \cite[Theorem 1]{AOB19}, which states that $v_\la$ is the unique solution of $\val(W_\la(z))=0$, and from the theory of matrix games \cite{SS50} which implies the existence of a square sub-matrix corresponding to the support an a pair of extreme optimal strategies. \hfill $\square$


\subsection{Proof of Theorem 1}
Consider an absorbing game $(g,q,w)\in \QQ^{m\times n}$. By Lemma \ref{lem}, 
the discounted value $v_\la$ is the solution of $\det(\dot{W}_\la(z))=0$ for some square sub-matrix $\dot{W}_\la(z)$ of size $r\leq \min(m,n)$. 
Assuming $\min(m,n)<3$, this leads to consider only sub-matrices of size $1$ and $2$. In the first case, this means the existence of $(i,j)$ so that $W_\la(z)_{ij}=0$ and for all $\la$ sufficiently small. Hence: 
$$v_\la= \frac{\la g_{ij} +(1-\la)q_{ij} w_{ij}}{\la + (1-\la)q_{ij}},\quad \forall \la\in (0,\la_0)\,.$$
 Taking $\la$ to $0$ gives then $v=\lim_{\la\to 0}v_\la=  w_{ij}\in \QQ$. \\
Second, suppose that $v_\la$ is the solution of a minor of size 2, say $\{1,2\}\times\{1,2\}$ w.l.o.g. (up to relabeling the actions). In this case, $v_\la$ is the solution of: 
$$z\in \RR,\quad \left| \begin{matrix}
W_\la(z)_{11} & W_\la(z)_{12} \\
W_\la(z)_{21} & W_\la(z)_{22} 
\end{matrix}\right|=0\,.$$
Taking the limit as $\la$ goes to $0$ gives then the following equation:
\begin{equation}\label{qw}q_{11}(w_{11}-z)q_{22}(w_{22}-z)=q_{12}(w_{12}-z)q_{21}(w_{21}-z)\,.\end{equation}
For $z$ to be an irrational solution, a necessary and sufficient condition is that the coefficient of $z^2$ is non-zero, and that both roots of degree-$2$ polynomial are irrational. Let us show that this is impossible to achieve for absorbing games with deterministic transitions, i.e., $q_{i,j}\in \{0,1\}$ for all $(i,j)$. Indeed, the coefficient for $z^2$ is $q_{11}q_{22}-q_{12}q_{21}$. To be different from $0$, either $q_{11}q_{22}=1$ and $q_{12}q_{21}=0$, or the converse. Replacing these values in equation \eqref{qw} gives then $(w_{11}-z)(w_{22}-z)=0$ or $(w_{12}-z)(w_{21}-z)=0$, respectively, which have only rational solutions.\hfill $\square$ 
\subsection{Proof of Theorem 2}
We provide two proofs: first, one that goes straight to the point but gives no intuition on how the example was found; second, one that is more constructive. 

\paragraph{First proof.} By playing a fixed stationary strategy $x=(x^1,x^2,x^3)$, Player $1$ can guarantee $\min_{j\in J} \lim_{\la\to 0} \ga_\la(x,j)$. Maximizing over $\De(I)$ gives then:
\[v\geq \max_{x=(x^1,x^2,x^3)\in \De(I)}\min \left\{x^1 + x^2 + 2x^3, \frac{x^1+2x^2}{x^1+x^2},\frac{2x^1+x^3}{x^1+x^3}\right\}\,.\]
Looking for an equalizing strategy (i.e., so that $\lim_{\la\to 0} \min_j \ga_\la^1(x,j)$ is independent of $j$) leads then to the following system of equations: 
\begin{equation}\label{eq1}
\left\{
\begin{aligned}
x^1 + x^2 + 2x^3 &= \frac{x^1+2x^2}{x^1+x^2}, \\
x^1 + x^2 + 2x^3 &= \frac{2x^1+x^3}{x^1+x^3}, \\
x^1 + x^2 + x^3 &= 1, \\
x^1, x^2, x^3 &\geq 0\,.
\end{aligned}
\right.
\end{equation}
This system has a unique solution\footnote{Solved with Wolfram Alpha, at https://www.wolframalpha.com (1 July 2023).} $x^1=\al$, $x^2=1-2\al-\al^2$, and $x^3=\al+\al^2$, where
$$\al = -1 + \frac{1}{3} \sqrt[3]{\frac{27}{2} - \frac{3 \sqrt{69}}{2}} + \frac{\sqrt[3]{\frac{1}{2}(9 + \sqrt{69})}}{\sqrt[3]9}\simeq 0.3247\,.$$
Using this equalizer stationary strategy, Player $1$ obtains $1+\al+\al^2$ no matter what Player 2 plays, hence 
$v\geq 1+\al+\al^2$. The symmetry of the game allows then to revert the roles of the players, and thus obtain %
\begin{equation}\label{expv}v=1+\al+\al^2= \frac{1}{3}\left(5 - 5 \sqrt[3]{\frac{ 2}{ 3 \sqrt {69} - 11}} + \sqrt[3]{\frac{1}{2}(3 \sqrt {69} - 11)}\right) \simeq 1.43\,.\end{equation}
Note that the value, but also each coordinate of a limit optimal stationary strategies, are algebraic of order $3$.\hfill $\square$  

\paragraph{Second proof.}
We proceed in 3 steps. \\
\noindent \textit{Step 1.} 
Consider the stationary strategy $x=(1/3,1/3,1/3)$. Then, 
$$v\geq \lim_{\la\to 0 } \min_{j} \ga_{\la}(x,j)= \lim_{\la\to 0 } \ga_{\la}(x,1)= \frac 4 3\,. $$
\noindent \textit{Step 2.} 
Let $(x_\la)$ be a family of optimal stationary strategies of Player 1 and let $x$ be an accumulation point. Then, $x^i>0$ for all $i$. 
Indeed, assume on the contrary that $x_{\la_n}$ tends to $x=(x^1,x^2,x^3)$ where $x^1=0$, $x^2=0$ or $x^3=0$. Consider the cases separately: for each of them, there exists a pure action of Players $2$ (namely, $j=3$, $j=2$, and $j=1$, respectively) so that $$\lim_{n\to +\infty } \ga^1_{\la_n}(x_{\la_n},j)\leq 1.$$
By the choice of $(x_{\la_n})$, this entails then
$$v=\lim_{n\to +\infty }v_{\la_n}=\lim_{n\to +\infty }\min_{j} \ga^1_{\la_n}(x_{\la_n},j)\leq1,$$ 
which contradicts the result of Step $1$. Hence, Player 1's optimal strategies have full support for all sufficiently small $\la$. \\ 
\noindent \textit{Step 3.} The $W_\la(z)$ matrix corresponding to our example is: 
\[
W_\la^1(z) =\begin{bmatrix}
1-z & 1-z & 2-z \\
1-z & 2-z & -\la z \\
2-z & -\la z  & 1-z
\end{bmatrix}\,.
\] 
Given that Player 1's optimal strategies have full support for all sufficiently small $\la$, Lemma \ref{lem} implies that  $P(\la,z):=\det(W_\la(z))$ is a polynomial of degree at least $1$ in $z$ satisfying $P(\la,v_\la)=0$ for all $\la$ sufficiently small. Taking $\la$ to $0$, one thus proves that $v$ is a root of the following polynomial
$$P_0(z):=P(0,z)=\det\begin{bmatrix}
1-z & 1-z & 2-z \\
1-z & 2-z & 0 \\
2-z & 0  & 1-z
\end{bmatrix}= z^3-5z^2+10z-7\,.$$ This polynomial has a unique real solution, given by equation \eqref{expv}, which is then equal to $v$. \\ 

\subsection{Proof of Theorem 3}
\paragraph{First proof.} Like in the first proof of Theorem 1, by playing a fixed stationary strategy $x=(x^1,x^2)$, Player $1$ can guarantee $\min_{j\in J} \lim_{\la\to 0} \ga_\la^1(x,j)$. Maximizing over $\De(I)$ gives then:
\[v\geq \max_{x=(x^1,x^2)\in \De(I)}\min \left\{\frac{k}{x^1+k x^2},x^1+k x^2\right\}\,.\]
Looking for a equalizer yields to the following system of equations:
\begin{equation}\label{eq2} 
\left\{
\begin{aligned}
\frac{k}{x^1+k x^2}&=x^1+k x^2,\\
x^1+x^2&=1,\\
x^1,x^2&\geq 0\,.
\end{aligned}
\right.
\end{equation}
which admits a unique solution $x_0=(\frac{\sqrt k -1}{k-1}, \frac{k - \sqrt k}{k-1})$. This strategy ensures $\sqrt{k}$ no matter what Player 2 plays, so that $v\geq \sqrt k$.
By symmetry, the same strategy is equalizer for Player $2$, so that $v\leq\sqrt k$, and hence $v=\sqrt k$.\hfill $\square$
\paragraph{Second proof.} Alternatively, one can start by proving that Player 1's optimal strategies has full support: otherwise it can only guarantee $1$, which is less than $\lim_{\la\to 0}\min_j \ga_\la(1/2,1/2; j)$ for all $k\geq 2$. Hence, by Lemma \ref{lem}, $v_\la$
solves
$$z\in \RR, \quad \det(W_\la(z))=\frac{1}{k}(k-z)^2-(1-z)^2=0\,.$$
The unique positive solution of this equation is $\sqrt{k}$, which concludes the proof. \hfill $\square$ 


\section{Discussion}
\begin{itemize}
    \item The irrationality of the limit value in this context comes from the non-linearity the limit payoffs, which appear for example in the equalizing equations \eqref{eq1} and \eqref{eq2}.
\item To any absorbing game with $m$ actions per player corresponds the parameterized (limit) matrix $A(z):=\lim_{\la \to 0} W_\la(z)$ given by 
$$A(z)_{ij}\in \{q_{ij}(w_{ij}-z), 0\},\quad \forall (i,j)\,.$$ Can we obtain any desired polynomial with rational coefficients as $\det(A(z))$ by varying $m$ and $q,w\in \RR^{m\times m}$? 
\end{itemize} 


\section{Acknowledgements}
The author is grateful to Eilon Solan, Guillaume Vigeral, and Krishnendu Chatterjee for valuable discussions that helped improve this paper. This work was supported by the French Agence Nationale de la Recherche (ANR) under reference ANR-21-CE40-0020 (CONVERGENCE project). 

\bibliographystyle{plain}
\bibliography{biblio}


\end{document}